\documentclass[11pt]{articlefederico}

\usepackage{amsthm}
\usepackage{amsmath}
\usepackage{amssymb}

\usepackage[left=1.5in,top=1.0in,right=1.5in]{geometry}

\usepackage{graphicx}

\newtheorem{theorem}[subsection]{Theorem}
\newtheorem{cor}[subsection]{Corollary}
\newtheorem{lemma}[subsection]{Lemma}
\newtheorem{prop}[subsection]{Proposition}

\theoremstyle{remark}
\newtheorem{exa}[subsection]{Example}

\theoremstyle{definition}
\newtheorem{deff}[subsection]{Definition}

\newcommand{\conv}{\mathrm{conv}}
\newcommand{\R}{\mathbb{R}}
\newcommand{\vol}{\mathrm{Vol\,}}
\newcommand{\x}{\mathbf{x}}
\newcommand{\bbeta}{\widetilde{\beta}}
\newcommand{\ggamma}{\widetilde{\gamma}}

\newcommand{\C}{\mathbb{C}}

\newcommand{\F}{\mathcal{F}}

\begin{document}
\title{\textsf{matroid polytopes and their volumes.}}
\author{\textsf{federico ardila\footnote{\textsf{San Francisco State University, San Francisco, CA, USA, federico@math.sfsu.edu.}}}
\qquad \textsf{carolina benedetti\footnote{\textsf{Universidad de Los Andes, Bogot\'a, Colombia, c.benedetti92@uniandes.edu.co.}}}
\qquad \textsf{jeffrey doker\footnote{\textsf{University of California, Berkeley, Berkeley, CA, USA, doker@math.berkeley.edu.
\newline
This research was partially supported by NSF grant DMS-0801075 (Ardila), the Proyecto Semilla of the Universidad de Los Andes (Benedetti), and a San Francisco State University Seed Funding Grant to support the SFSU-Colombia Combinatorics Initiative.}}}}

\date{}

\maketitle

\begin{abstract}
We express the matroid polytope $P_M$ of a matroid $M$ as a signed Minkowski sum of simplices, and obtain a formula for the volume of $P_M$. This gives a combinatorial expression for the degree of an arbitrary torus orbit closure in the Grassmannian $Gr_{k,n}$. We then derive analogous results for the independent set polytope and the associated flag matroid polytope of $M$. Our proofs are based on a natural extension of Postnikov's theory of generalized permutohedra.
\end{abstract}

\section{\textsf{introduction.}} \label{sec:introduction}

The theory of matroids can be approached from many different points of view; a matroid can be defined as a simplicial complex of independent sets, a lattice of flats, a closure relation, etc. A relatively new point of view is the study of matroid polytopes, which in some sense are the natural combinatorial incarnations of matroids in algebraic geometry and optimization. Our paper is a contribution in this direction.

We begin with the observation that matroid polytopes are members of the family of \emph{generalized permutohedra} \cite{Po}. With some modifications of Postnikov's beautiful theory, we express the matroid polytope $P_M$  as a signed Minkowski sum of simplices, and use that to give a formula for its volume $\vol(P_M)$. This is done in Theorems \ref{th:sum} and \ref{th:volume}. Our answers are expressed in terms of the beta invariants of the contractions of $M$. 

Formulas for $\vol(P_M)$ were given in very special cases by Stanley \cite{St} and Lam and Postnikov \cite{LP}, and a polynomial-time algorithm for finding $\vol(P_M)$ was constructed by de Loera et. al. \cite{de}. One motivation for this computation is the following. The closure of the torus orbit of a point $p$ in the Grassmannian $Gr_{k,n}$ is a toric variety $X_p$, whose degree is the volume of the matroid polytope $P_{M_p}$ associated to $p$. Our formula allows us to compute the degree of $X_p$ combinatorially.

One can naturally associate two other polytopes to a matroid $M$: its independent set polytope and its associated flag matroid polytope. By a further extension of Postnikov's theory, we also write these polytopes as signed Minkowski sums of simplices and give formulas for their volumes. This is the content of Sections \ref{sec:indep} and \ref{sec:flagmatroids}.

Throughout the chapter we assume familiarity with the basic concepts of matroid theory; for further information we refer the reader to \cite{Ox}.

\section{\textsf{matroid polytopes are generalized permutohedra}} \label{sec:permutohedra}

The \emph{permutohedron} $P_n$ is a polytope in $\mathbb{R}^n$ whose vertices consist of all permutations of the entries of the vector $(1,2,\dots,n)$. A \emph{generalized permutohedron} is a deformation of the permutohedron, obtained by moving the vertices of $P_n$ in such a way that all edge directions and orientations are preserved (and some may possibly be shrunken down to a single point) ~\cite{PRW}.

Every generalized permutohedron can be written in the following form:
\[
P_n(\{z_I\}) = \left\{(t_1, \ldots, t_n) \in \R^n : \sum_{i=1}^n t_i =
z_{[n]}, \sum_{i \in I} t_i \geq z_I \textrm{ for all } I
\subseteq [n]\large\right\}
\]
where $z_I$ is a real number for each $I \subseteq [n]:=\{1, \ldots, n\}$, and
$z_{\emptyset}=0$. Different choices of $z_I$ can give the same
generalized permutohedron: if one of the inequalities does not
define a face of $P_n(\{z_I\})$, then we can decrease the value
of the corresponding $z_I$ without altering the polytope. When we
write $P_n(\{z_I\})$, we will always assume that the $z_I$s are all
chosen maximally; \emph{i.e.}, that all the defining inequalities
are tight.

Though every generalized permutohedron has a $z_I$ parameterization, not every list of $z_I$ parameters corresponds to a generalized permutohedron. 
The following characterization was announced by Morton et. al. \cite[Theorem 17]{Morton} and Postnikov \cite{Po2}. For a complete proof, see \cite{AA}.

\begin{theorem}
\label{thm: Morton}
A set of parameters $\{z_I\}$ defines a generalized permutohedron $P_n(\{z_I\})$ if and only if the $z_I$ satisfy the supermodular inequalities
\[
z_I + z_J \leq z_{I\cup J} + z_{I\cap J}
\]
for all $I,J\subseteq [n]$.
\end{theorem}



The \emph{Minkowski sum} of two polytopes $P$ and $Q$ in $\R^n$
is defined to be $P+Q = \{p+q \, : \, p \in P, q \in Q\}$.
We say that the \emph{Minkowski difference} of $P$ and $Q$ is $P-Q=R$ if
$P=Q+R$.\footnote{We will only consider Minkowski differences
$P-Q$ such that $Q$ is a Minkowski
summand of $P$. More generally, the Minkowski difference of
two arbitrary polytopes $P$ and $Q$ in
$\R^n$  is defined to be $P-Q = \{r\in \R^n \, | \, r+Q \subseteq
P\}$ \cite{Po}. It is easy to check that $(Q+R)-Q=R$, so the
two definitions agree in the cases that interest us. In this
paper, a signed Minkowski sum equality such as
$P-Q+R-S=T$ should be interpreted as $P+R=Q+S+T$.}
The following lemma shows that generalized permutohedra
behave nicely with respect to Minkowski sums.

\begin{lemma}\label{lemma:add}
If $P_n(\{z_I\})$ and $P_n(\{z'_I\})$ are generalized permutohedra then their Minkowski sum is also a generalized permutohedron and $P_n(\{z_I\}) + P_n(\{z'_I\}) = P_n(\{z_I+z'_I\})$.
\end{lemma}

\begin{proof}
The polytopes $P_n(\{z_I\})$ and $P_n(\{z'_I\})$ are deformations of $P_n$, and therefore by \cite[Theorem 17]{Morton} they are each a Minkowski summand of a dilate of $P_n$. Thus $P_n(\{z_I\}) + P_n(\{z'_I\})$ must also be a summand of a dilate of $P_n$, which implies, again by \cite[Theorem 17]{Morton}, that this polytope too is a deformation of $P_n$ and can thus be defined by hyperplane parameters $z_I$. That the values of these parameters are $z_I+z'_I$ follows from the observation that, if a linear functional $w$
takes maximum values $a$ and $b$ on (faces $A$ and $B$ of)
polytopes $P$ and $Q$ respectively, then it takes maximum value
$a+b$ on (the face $A+B$ of) their Minkowski sum. 
\end{proof}

Let $\Delta$ be the standard unit $(n-1)$-simplex
\begin{eqnarray*}
\Delta&=&\{(t_1, \ldots, t_n) \in \R^n\, : \, \sum_{i=1}^n t_i = 1, t_i \geq 0 \textrm{ for all } 1 \leq i \leq n\}\\
&=& \conv\{e_1, \ldots, e_n\},
\end{eqnarray*}
where $e_i=(0,,\ldots, 0, 1, 0, \ldots, 0)$ with a $1$ in its
$i$th coordinate. As $J$ ranges over the subsets of $[n]$, let
$\Delta_J$ be the face of the simplex $\Delta$ defined by
\[
\Delta_J = \conv\{e_i \, : \, i \in J\} = P_n(\{z(J)_I\})
\]
where $z(J)_I=1$ if $I \supseteq J$ and $z(J)_I=0$ otherwise. Lemma \ref{lemma:add} gives the following proposition.

\begin{prop}\label{prop:sum.facets}\cite[Proposition 6.3]{Po}
For any $y_I \geq 0$, the Minkowski sum $\sum y_I \Delta_I$ of dilations of faces of the standard $(n-1)$-simplex is a
generalized permutohedron. We can write
\[
\sum_{A \subseteq E} y_I \Delta_I = P_n(\{z_I\}),
\]
where $z_I = \sum_{J \subseteq I} y_J$ for each $I \subseteq [n]$.
\end{prop}

We can extend this to encompass signed Minkowski sums as well.

\begin{prop}\label{prop:facets.sum}
Every generalized permutohedron $P_n(\{z_I\})$ can be written uniquely as a signed Minkowski sum of simplices, as
\[
P_n(\{z_I\}) = \sum_{I \subseteq [n]} y_I \Delta_I
\]
where $y_I = \sum_{J \subseteq I} (-1)^{|I|-|J|} z_J$ for each $I \subseteq [n]$.
\end{prop}

\begin{proof}
First we need to separate the righthand side into its positive and negative parts. By Proposition \ref{prop:sum.facets},
\[
\sum_{I \subseteq [n] \, : y_I < 0} (-y_I) \Delta_I =
P_n(\{z^-_I\}) \textrm{ and } \sum_{I \subseteq [n] \, : y_I \geq
0} y_I \Delta_I = P_n(\{z^+_I\})
\]
where $z^-_I= \sum_{J \subseteq I \, : y_J < 0} (-y_J)$ and $z^+_I
= \sum_{J \subseteq I \, : y_J \geq 0} y_J$. Now 
$z_I +z^-_I = z^+_I$ gives
\[
P_n(\{z_I\}) + \sum_{I \subseteq [n] \, : y_I < 0} (-y_I) \Delta_I
= \sum_{I \subseteq [n] \, : y_I \geq 0} y_I \Delta_I,
\]
as desired. Uniqueness is clear.
\end{proof}

Let $M$ be a matroid of rank $r$ on the set $E$. The
\emph{matroid polytope} of $M$ is the polytope $P_M$ in $\R^E$
whose vertices are the indicator vectors of the bases of $M$.
The known description of the polytope $P_M$ by inequalities makes it apparent
that it is a generalized permutohedron:
\begin{prop}\label{prop:facets.matroid}\cite{We}
The matroid polytope of a matroid $M$ on $E$ with rank function $r$ is
$P_M=P_E(\{r-r(E-I)\}_{I \subseteq E}).$
\end{prop}
\begin{proof}
The inequality description for $P_M$ is:
\[
P_M = \{\x \in \R^E \, : \, \sum_{i \in E} x_i = r, \sum_{i \in A}
x_i \leq r(A) \textrm{ for all } A \subseteq E\}.
\]
It remains to remark that the inequality $\sum_{i \in A} x_i \leq
r(A)$ is tight, and may be rewritten as $\sum_{i \in E-A} x_i \geq
r-r(A)$, and to invoke the submodularity of the rank function of a matroid.
\end{proof}

The \emph{beta invariant} \cite{Cr} of $M$ is a non-negative
integer given by
\[
\beta(M) = (-1)^{r(M)} \sum_{X \subseteq E} (-1)^{|X|} r(X)
\]
which stores significant information about $M$; for example, $\beta(M)=0$ if and only if $M$ is disconnected and $\beta(M)=1$ if and only if $M$ is series-parallel. If
\[
T_M(x,y) = \sum_{A \subseteq E} (x-1)^{r(E)-r(A)} (y-1)^{|A|-r(A)} = \sum_{i,j} b_{ij}x^iy^j
\]
is the \emph{Tutte polynomial} \cite{Wh} of $M$, then $\beta(M)=b_{10}=b_{01}$ for $|E|\geq 2$.

Our next results are more elegantly stated in terms of the
\emph{signed beta invariant} of $M$, which we define to be
\[
\bbeta(M) = (-1)^{r(M)+1} \beta(M).
\]

\begin{theorem}\label{th:sum}
Let $M$ be a matroid of rank $r$ on $E$ and let $P_M$ be its matroid
polytope. Then
\begin{equation}\label{eq:sum}
P_M = \sum_{A \subseteq E}  \bbeta(M/A)\,
\Delta_{E-A}.
\end{equation}
\end{theorem}

\begin{proof}
By Propositions \ref{prop:facets.sum} and
\ref{prop:facets.matroid}, $P_M = \sum_{I \subseteq E} y_I
\Delta_I$ where
\begin{eqnarray*}
y_I &=&  \sum_{J \subseteq I} (-1)^{|I|-|J|} (r - r(E-J))
=  -\sum_{J \subseteq I} (-1)^{|I|-|J|} r(E-J)\\
&=&  -\sum_{E-J \supseteq E-I} (-1)^{|E-J|-|E-I|} (r(E-J)-r(E-I))\\
&=&  -\sum_{X \subseteq I} (-1)^{|X|} (r(E-I\cup X)-r(E-I))\\
&=& -\sum_{X \subseteq I} (-1)^{|X|} r_{M/(E-I)}(X)= \bbeta(M/(E-I))
\end{eqnarray*}
as desired.
\end{proof}

\begin{figure}[h]
\centering
\includegraphics[scale=.68]{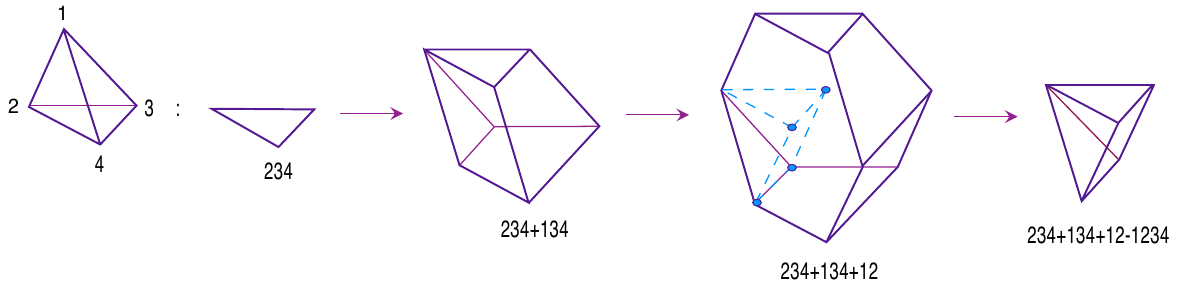} 
\caption{A matroid polytope as a signed Minkowski sum of simplices.\label{fig1}}
\end{figure}

\begin{exa}
Let $M$ be the matroid on $E=[4]$ with bases $\{12,13,14,23,24\}$; its matroid polytope is a square pyramid. Theorem \ref{th:sum} gives $P_M=\Delta_{234}+\Delta_{134}+\Delta_{12}-\Delta_{1234}$, as illustrated in Figure \ref{fig1}. The dotted lines in the polytope $\Delta_{234}+\Delta_{134}+\Delta_{12}$ are an aid to visualize the Minkowski difference.

\end{exa}

One way
of visualizing the Minkowski sum of two polytopes $P$ and
$Q$ is by grabbing a vertex $v$ of $Q$ and then using it to ``slide''
$Q$ around in space, making sure that $v$ never leaves $P$.  The
region that $Q$ sweeps along the way is $P+Q$.
%
Similarly, the Minkowski
difference $P-R$ can be visualized by picking a vertex $v$ of $R$ and
then ``sliding'' $R$ around in space, this time making sure that no
point in $R$ ever leaves $P$. The region that $v$ sweeps along the
way is $P-R$.
This may be helpful in understanding Figure \ref{fig1}. 

Some remarks about Theorem \ref{th:sum} are in order.

\begin{itemize}
\item  Generally most terms in the sum of Theorem \ref{th:sum} are zero. The nonzero terms correspond to the \emph{coconnected flats} $A$, which we define to be the sets $A$ such that $M/A$ is connected. 
These are indeed flats, since contracting by them must
produce a loopless matroid. 


\item
A matroid and its dual have congruent matroid polytopes, and
Theorem \ref{th:sum} gives different
formulas for them. 
For example $P_{U_{1,3}} = \Delta_{123}$ while 
$P_{U_{2,3}} = \Delta_{12} + \Delta_{23} + \Delta_{13} - \Delta_{123}$.

%




\item
The study of the subdivisions of a matroid polytope into smaller matroid polytopes, originally considered by Lafforgue \cite{La}, has recently received significant attention \cite{AFR, BJR, De, Sp04}. Speyer conjectured  \cite{Sp04} that the subdivisions consisting of series-parallel matroids have the largest number of faces in each dimension and proved this \cite{Sp07} for a large and important family of subdivisions: those that arise from a tropical linear space. 
The important role played by series-parallel matroids is still somewhat mysterious.
Theorem \ref{th:sum} characterizes series-parallel matroids as
those whose matroid polytope has no repeated Minkowski summands. It would be interesting to connect this characterization to matroid subdivisions; this may require extending the theory of mixed subdivisions to signed Minkowski sums.

\item
Theorem \ref{th:sum} provides a geometric interpretation for the
beta invariant of a matroid $M$ in terms of the
matroid polytope $P_M$. In Section \ref{sec:flagmatroids}
we see how to extend this to certain families of Coxeter matroids.
This is a promising point of view
towards the notable open problem \cite[Problem
6.16.6]{BGW} of defining useful enumerative invariants
of a Coxeter matroid.
\end{itemize}

\section{\textsf{the volume of a matroid polytope}}\label{sec:volume}

Our next goal is to present an explicit combinatorial formula for the
volume of an arbitrary matroid polytope. Formulas have been
given for very special families of matroids by Stanley
\cite{St} and Lam and Postnikov \cite{LP}. Additionally, a
polynomial time algorithm for computing the volume of an arbitrary
matroid polytope was recently discovered by de Loera et. al. \cite{de}. 
Let us say some words about the motivation for this question.

Consider the Grassmannian manifold $Gr_{k,n}$ of $k$-dimensional
subspaces in $\C^n$. Such a subspace can be represented as the
rowspace of a $k \times n$ matrix $A$ of rank $k$, modulo the
left action of $GL_k$ which does not change the row space.
The ${n \choose k}$ maximal minors of this matrix are the
\emph{Pl\"ucker coordinates} of the subspace, and they give
an embedding of $Gr_{k,n}$ as a projective algebraic variety in $\C\mathbb{P}^{{n \choose k} -1}$.

Each point $p$ in $Gr_{k,n}$ gives rise to a matroid $M_p$ whose
bases are the $k$-subsets of $n$ where the Pl\"ucker coordinate
of $p$ is not zero. 
Gelfand, Goresky, MacPherson, and Serganova \cite{GGMS}
first considered the stratification of $Gr_{k,n}$ into
\emph{matroid strata}, which consist of the points
corresponding to a fixed matroid. 

The torus $\mathbb{T} = (\C^*)^n$ acts on $\C^n$ by
$(t_1, \ldots, t_n) \cdot (x_1, \ldots, x_n) = (t_1x_1, \ldots, t_nx_n)$
for $t_i \neq 0$; this action extends to an action of $\mathbb{T}$ on
$Gr_{k,n}$. For a point $p \in Gr_{k,n}$, the closure of the torus
orbit $X_p = \overline{\mathbb{T} \cdot p}$ is a toric variety which
only depends on the matroid $M_p$ of $p$, and the polytope corresponding
to $X_p$ under the moment map is the matroid polytope of $M_p$ \cite{GGMS}. Under these
circumstances it is known \cite{Fu} that the volume of the matroid
polytope $M_p$ equals the degree of the toric variety $X_p$ as a
projective subvariety of $\C \mathbb{P}^{{n \choose k} -1}$:
\[
\vol P_{M_p} = \textrm{deg}\, X_p.
\]
Therefore, by finding the volume of an arbitrary matroid polytope,
one obtains a formula for the degree of the toric varieties arising
from arbitrary torus orbits in the Grassmannian.

\medskip

To prove our formula for the volume of a matroid polytope, we
first recall the notion of the \emph{mixed
volume} $\vol(P_1, \ldots, P_n)$ of $n$ (possibly repeated)
polytopes $P_1, \ldots, P_n$ in $\R^n$. All volumes in this section
are normalized with respect to the lattice generated by $e_1-e_2,
\ldots, e_{n-1}-e_n$ where our polytopes live; so the standard simplex
$\Delta$ has volume $1/(n-1)!$.

\begin{prop}\cite{Mc}\label{prop:mixed}
\label{prop:mixed.vol} Let $n$ be a fixed positive integer. There
exists a unique function $\vol(P_1, \ldots, P_n)$ defined on
$n$-tuples of polytopes in $\R^n$, called the \emph{mixed volume}
of $P_1, \ldots, P_n$, such that, for any collection of polytopes
$Q_1, \ldots, Q_m$ in $\R^n$ and any nonnegative real numbers
$y_1, \ldots, y_m$, the volume of the Minkowski sum $y_1Q_1 +
\cdots + y_mQ_m$ is the polynomial in $y_1, \ldots, y_m$ given by
\[
\vol(y_1Q_1+ \cdots + y_mQ_m) = \sum_{i_1, \ldots, i_n}
\vol(Q_{i_1}, \ldots, Q_{i_n}) y_{i_1}\cdots y_{i_n}
\]
where the sum is over all ordered $n$-tuples $(i_1, \ldots, i_n)$
with $1 \leq i_r \leq m$.
\end{prop}

We now show that the formula of Proposition
\ref{prop:mixed.vol} still holds if some of the $y_i$s are
negative as long as the expression  $y_1Q_1+ \cdots + y_mQ_m$ still makes
sense.

\begin{prop}\label{prop:gen.mixed.vol}
If $P=y_1Q_1+ \cdots + y_mQ_m$ is a \emph{signed} Minkowski sum of
polytopes in $\R^n$, then
\[
\vol(y_1Q_1+ \cdots + y_mQ_m) = \sum_{i_1, \ldots, i_n}
\vol(Q_{i_1}, \ldots, Q_{i_n}) y_{i_1}\cdots y_{i_n}
\]
where the sum is over all ordered $n$-tuples $(i_1, \ldots, i_n)$
with $1 \leq i_r \leq m$.
\end{prop}

\begin{proof}
We first show that
\begin{equation}\label{bla1}
\vol (A-B) = \sum_{k=0}^n (-1)^k {n \choose k} \vol(A, \ldots, A, B, \ldots, B)
\end{equation}
when $B$ is a Minkowski summand of $A$ in $\R^n$. Let $A-B=C$. By
Proposition \ref{prop:mixed.vol}, for $t \geq 0$ we have that
\[
\vol(C+tB) = \sum_{k=0}^n {n \choose k} \vol(C, \ldots, C, B \ldots, B) t^k =:
f(t)
\]
and we are interested in computing $\vol(C)=f(0)$. Invoking
Proposition \ref{prop:mixed.vol} again, for $t \geq 0$ we have that
\begin{equation}\label{bla2}
\vol(A+tB) = \sum_{k=0}^n {n \choose k} \vol(A, \ldots, A, B, \ldots, B) t^k =:
g(t).
\end{equation}
But $A+tB = C+(t+1)B$ and therefore $g(t)=f(t+1)$ for all $t \geq 0$.
Therefore $g(t)=f(t+1)$ as polynomials, and
$\vol C = f(0) = g(-1)$. Plugging into (\ref{bla2}) gives the desired result.

Having established (\ref{bla1}), separate the given Minkowski sum for $P$
into its positive and negative parts as $P=Q-R$, where $Q=x_1Q_1 +
\cdots + x_rQ_r$ and $R = y_1R_1 +\cdots+y_sR_s$ with $x_i, y_i
\geq 0$. For positive $t$ we can write $Q+tR = \sum x_i Q_i + \sum
ty_j R_j$, which gives two formulas for $\vol(Q+tR)$.
\begin{eqnarray*}
\vol(Q+tR) &=& \sum_{k=0}^n {n \choose k} \vol(Q, \ldots, Q, R, \ldots, R)t^k\\
&=& \sum_{\stackrel{1 \leq i_a \leq r}{1 \leq j_b \leq s}}
 \vol(Q_{i_1}, \ldots, Q_{i_{n-k}},R_{j_1}, \ldots,
R_{j_k}) x_{i_1} \cdots x_{i_{n-k}} y_{j_1} \cdots y_{j_k} t^k
\end{eqnarray*}
The last two expressions must be equal as polynomials.
A priori, we cannot plug $t=-1$ into this equation; but instead, we can use
the formula for $\vol(Q-R)$ from (\ref{bla1}), and then plug in
coefficient by coefficient. That gives the desired result.
\end{proof}

\begin{theorem}\label{th:volume}
If a connected matroid $M$ has $n$ elements, then the volume of the matroid polytope $P_M$ is
\[
\vol P_M = \frac{1}{(n-1)!} \sum_{(J_1,
\ldots, J_{n-1})} \bbeta(M/J_1) \bbeta(M/J_2)
\cdots \bbeta(M/J_{n-1}),
\]
summing over the ordered collections of sets
$J_1, \ldots, J_{n-1} \subseteq [n]$ such that, for any distinct $i_1,
\ldots, i_k$, $|J_{i_1} \cap \cdots \cap J_{i_k}| < n - k$.
\end{theorem}

\begin{proof}
Postnikov \cite[Corollary 9.4]{Po} gave a formula for the
volume of a (positive) Minkowski sum of simplices. We would
like to apply his formula to the signed Minkowski sum in
Theorem \ref{th:sum}, and  Proposition \ref{prop:gen.mixed.vol} makes this possible.
\end{proof}

There is an alternative characterization of the
tuples $(J_1, \ldots, J_{n-1})$ considered in the sum above. They
are the tuples such that, for each $1 \leq k \leq n$, the
collection $([n]-J_1, \ldots, [n]-J_{n-1})$ has a system of
distinct representatives avoiding $k$; that is, there exist $a_1 \in
[n] - J_1, \ldots, a_{n-1} \in [n] - J_{n-1}$ with $a_i \neq a_j$ for $i \neq
j$ and $a_i \neq k$ for all $i$. Postnikov refers to this as the
\emph{dragon marriage condition}; see \cite{Po} for an explanation
of the terminology.

As in Theorem \ref{th:sum}, most of the terms in the sum of
Theorem \ref{th:volume} vanish. The nonzero terms are those
such that each $J_i$ is a coconnected flat. Furthermore,
since $P_M$ and $P_{M^*}$ are congruent, we are free to
apply Theorem \ref{th:volume} to the one giving a simpler expression.

\begin{exa}
Suppose we wish to compute the volume of $P_{U_{2,3}}$ using
Theorem \ref{th:volume}. The expression $P_{U_{1,3}} = \Delta_{123}$
is simpler than the one for $P_{U_{2,3}}$. So we can obtain $\vol P_{(U_{1,3})^*} = \vol P_{U_{1,3}} = \frac12\bbeta(M)^2 = \frac12$.
\end{exa}


In Theorem \ref{th:volume}, the hypothesis that $M$ is connected
is needed to guarantee that the matroid polytope $P_M$ has
dimension $n-1$. More generally, if we have
$M=M_1\oplus \cdots \oplus M_k$ then $P_{M} = P_{M_1} \times \cdots \times P_{M_k}$
so the ($(n-k)$-dimensional) volume of $P_M$ is $\vol P_M = \vol P_{M_1} \cdots \vol P_{M_k}$.


\section{\textsf{independent set polytopes}}\label{sec:indep}

In this section we show that our analysis of matroid polytopes can
be carried out similarly for the \emph{independent set
polytope} $I_M$ of a matroid $M$, which is the convex hull of 
the indicator vectors of the
independent sets of $M$. The inequality description of $I_M$ is
known to be:
\begin{equation}\label{IM}
I_M = \{(x_1, \ldots, x_n) \in \R^n \, : \, x_i \geq 0 \textrm{ for } i\in [n],
\sum_{i \in A} x_i \leq r(A) \textrm{ for all } A \subseteq E\}.
\end{equation}

This realization of the independent set polytope of a matroid is not a generalized
permutohedron. Instead, it is a \emph{Q-polytope}. The class of $Q$-polytopes are the deformations of the simple polytope $Q_n$ whose vertices are formed by all distinct permutations of entries of the vectors $(1,\dots,n), (0,2,\dots,n), \dots, (0,\dots,0,n)$, and $(0,\dots,0)$. After translation, every $Q$-polytope can be expressed in the form
\begin{equation}\label{Q ineqs}
Q_n(\{z_J\}) = \left\{(t_1, \ldots, t_n) \in \R^n : t_i \geq 0 \textrm{
for all } i \in [n], \sum_{i \in J} t_i \leq z_J \textrm{ for all
} J \subseteq [n]\right\}
\end{equation}
where $z_J$ is a non-negative real number for each $J \subseteq [n]$. Analogously to generalized permutohedra, the parameters $z_J$ which describe a $Q$-polytope $Q_n(\{z_J\})$ are those which satisfy a supermodular inequality. However, since our inequalities are reversed in the hyperplane description of $Q_n$, we need to reverse the inequality on our defining supermodular criterion to obtain the submodular criterion.

\begin{prop}
The polytope $Q_n(\{z_J\})$ is a $Q$-polytope if and only if
\[
z_I + z_J \geq z_{I\cup J} + z_{I\cap J}
\]
for all $I,J\subseteq [n]$.
\end{prop}
\begin{proof}
This follows  from Theorem \ref{thm: Morton}. We will describe a deformation-preserving bijection from generalized permutohedra to $Q$-polytopes where supermodularity of $z_I$ parameters of $P_n(\{z_I\})$ corresponds to submodularity of $z_J$ parameters of ${Q}_n(\{z_J\})$. Given a generalized permutohedron $P_n(\{z_I\})\subset \R^n$, define $P_{n+1}(\{z'_I\})$ to be the generalized permutohedron in $\R^{n+1}$ defined by $z'_I = 0$ and $z'_{I\cup\{n+1\}} = z_I$ for all $I\subseteq[n]$. Now define ${Q}_n(\{z''_I\})$ to be the projection of $P_{n+1}(\{z'_I\})$ into $\R^n$, by removal of the last coordinate. This invertible process sends $P_n(\{z_I\})$ to ${Q}_n(\{z''_I\})$ where $z''_I = z_{[n]} - z_{[n]\setminus I}$. Moreover, it sends the permutohedron $P_n$ to $Q_n$ and deformations of $P_n$ to deformations of $Q_n$. A polytope $P$ has inequality description $P_n(\{z_I\})$ satisfying the supermodular inequalities if and only if $P$ is a deformation of the permutohedron $P_n$. By the map described above this occurs if and only if the corresponding polytope $Q = {Q}_n(\{z''_I\})$ is a deformation of $Q_n$. By the construction of the $z''_I$, supermodularity of the $z_I$ is equivalent to submodularity of the $z''_I$. Thus we have a deformation of $Q_n$ if and only if the corresponding $z''_I$ parameters are submodular.
\end{proof}

We can also express these polytopes as  signed
Minkowski sums of simplices, though the simplices we use are not
the $\Delta_J$s, but those of the form
\begin{eqnarray*}
D_J &=& \conv\{0, \,\, e_i \, : \, i \in J\}  \\
&=& Q_n(\{d(J)_I\})
\end{eqnarray*}
where $d(J)_I = 0$ if $I \cap J = \emptyset$ and $d(J)_I = 1$ otherwise.

The following lemmas on Q-polytopes are proved in a way
analogous to the corresponding lemmas for generalized
permutahedra, as was done in Section \ref{sec:permutohedra}.

\begin{lemma}
If $Q_n(\{z_J\})$ and $Q_n(\{z'_J\})$ are $Q$-polytopes, then so is their Minkowski sum, and $Q_n(\{z_J\})+Q_n(\{z_J'\})=Q_n(\{z_J+z_J'\})$
\end{lemma}


\begin{prop} For any $y_I \geq 0$ we have
\[\sum_{I \subseteq [n]} y_ID_I = Q_n(\{z_J\})\]
where $z_J=\sum_{I: I \cap J \neq \emptyset}y_I$.
\end{prop}
%



\begin{prop}\label{prop:Qn}
Every Q-polytope $Q_n(\{z_J\})$ can be written uniquely\footnote{assuming $y_{\emptyset}=0$} as a
signed Minkowski sum of $D_I$s as
\[Q_n(\{z_J\})=\sum_{I\subseteq [n]}y_ID_I, \]
where
\[y_J=-\sum_{I \subseteq J}(-1)^{|J|-|I|}z_{[n]-I}.\]
\end{prop}

\begin{proof}
We need to invert the relation between the $y_I$s and the
$z_J$s given by $z_J=\sum_{I: I \cap J \neq \emptyset}y_I$. We rewrite this relation as
\[
z_{[n]} - z_J =  \sum_{I \subseteq [n]-J}y_I
\]
and apply inclusion-exclusion.
As in Section \ref{sec:permutohedra}, we first do this in
the case $y_I \geq 0$ and then extend it to arbitrary Q-polytopes.
\end{proof}

\begin{theorem}\label{th:I_M}
Let $M$ be a matroid of rank $r$ on $E$ and let $I_M$
be its independent set polytope. Then
\[
I_M = \sum_{A \subseteq E} \bbeta(M/A) \, D_{E-A}
\]
where $\bbeta$ denotes the signed beta invariant.
\end{theorem}

\begin{proof}
This follows from Proposition \ref{prop:Qn} and a
computation almost identical to the one in the proof of Theorem \ref{th:sum}.
\end{proof}

The great similarity between Theorems \ref{th:sum}
and \ref{th:I_M} is not surprising, since $P_M$ is the
facet of $I_M$ which maximizes the linear function
$\sum_{i \in E} x_i$, and $\Delta_I$ is the facet of
$D_I$ in that direction as well. In fact we could have
first proved Theorem \ref{th:I_M} and then obtained
Theorem \ref{th:sum} as a corollary.

\begin{theorem}\label{vol:I_M}
If a connected matroid $M$ has $n$ elements, then the volume of the independent set polytope $I_M$ is
$$
\vol I_M=\frac{1}{n!} \sum_{(J_1, \ldots, J_{n})}
\bbeta(M/J_1) \bbeta(M/J_2) \cdots \bbeta(M/J_n)
$$
where the sum is over all $n-$tuples $(J_1,\dots,J_n)$ of subsets
of $[n]$ such that, for any distinct $i_1, \ldots, i_k$, we have $|J_{i_1} \cap \cdots \cap J_{i_k}| \leq n-k$.
\end{theorem}

Notice that by Hall's marriage theorem, the condition on
the $J_i$s is equivalent to requiring that $(E-J_1, \ldots, E-J_n)$
has a system of distinct representatives (SDR); that is,
there are $a_1 \in E-J_1, \ldots, a_n \in E-J_n$ with $a_i \neq a_j$ for $i \neq j$.

\begin{proof}
By Theorem \ref{th:I_M} and Proposition \ref{prop:mixed}
it suffices to compute the mixed volume $\vol(D_{A_1},\dots ,D_{A_n})$
for each $n$-tuple $(A_1, \ldots, A_n)$ of subsets of $[n]$.
Bernstein's theorem \cite{Stu} tells us that
$\vol(D_{A_1},\dots ,D_{A_n})$ is the number of isolated
solutions in $(\mathbb{C}-\{0\})^{n}$ of the system of equations:
\begin{align*}
\beta_{1,0}+\beta_{1,1}t_1+\beta_{1,2}t_2+\cdots +\beta_{1,n}t_n &=
0\\ \beta_{2,0}+\beta_{2,1}t_1+\beta_{2,2}t_2+\cdots +\beta_{2,n}t_n
&= 0\\ & \vdots \\
\beta_{n,0}+\beta_{n,1}t+\beta_{n,2}t_2+\cdots +\beta_{n,n}t_n &=
0\\
\end{align*}
where $\beta_{i,0}$ and $\beta_{i,j}$ are generic complex
numbers when $j \in A_i$, and $\beta_{i,j} = 0$ if $j\notin A_i$.

This system of linear equations will have one solution
if it is non-singular and no solutions otherwise. 
Because the $\beta_{i,0}$ are generic, such
a solution will be non-zero if it exists. The system is
non-singular when the determinant is non-zero, and by genericity
that happens when $(A_1, \ldots, A_n)$ has an SDR. We conclude that $\vol(D_{E-J_1},\dots ,D_{E-J_n})$ is $1$ if  $(E-J_1, \ldots, E-J_n)$ has an SDR and $0$ otherwise, and the result follows.
\end{proof}

Let us illustrate Theorem \ref{vol:I_M} with an example. 

\begin{figure}[h]
\centering
\includegraphics[scale=.4]{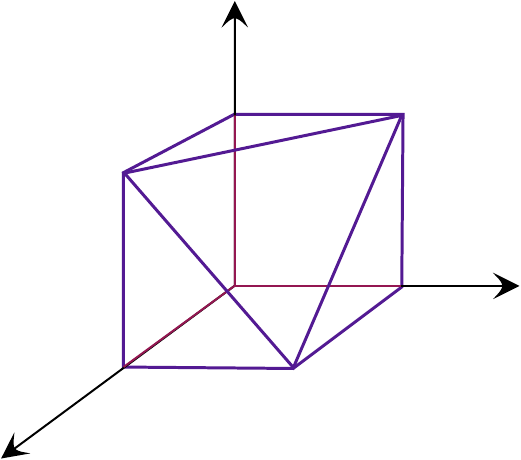} 
\caption{The independent set polytope of $U_{2,3}$.\label{fig2}}

\end{figure}

\begin{exa}
The independent set polytope of the uniform matroid $U_{2,3}$ is shown in Figure \ref{fig2}. We have $I_M = D_{12}+D_{23}+D_{13}-D_{123}$.
Theorem \ref{vol:I_M} should confirm that its volume is $\frac56$; let us
carry out that computation. 

The coconnected flats of $M$ are $1,2,3$ and $\emptyset$ and their
complements are $\{23,13,12,123\}$. We need to consider the triples of
coconnected flats whose complements contain an SDR. 
Each one of the 24 triples of the form $(a,b,c)$, where $a,b,c\in [3]$ are not all equal,
contributes a summand equal to 1. 
The 27 permutations of triples of the form $(a,b,\emptyset)$, contribute a
 $-1$ each.
The 9 permutations of triples of the form $(a,\emptyset, \emptyset)$ contribute
a 1 each.
The triple $(\emptyset,\emptyset, \emptyset)$ contributes a $-1$.
The volume of $I_M$ is then $\frac16(24-27+9-1) =\frac {5}{6}$.
\end{exa}

\section{\textsf{truncation flag matroids}}\label{sec:flagmatroids}

We will soon see that any flag matroid polytope can also be
written as a signed Minkowski sum of simplices $\Delta_I$.
We now focus on the particularly nice family of \emph{truncation flag matroids}, introduced by Borovik, Gelfand, Vince, and White \cite{BGVW}, where we obtain an explicit formula for this sum.

The \emph{strong order} on matroids is defined by saying that two matroids $M$ and $N$ 
on the same ground set $E$, having respective ranks $r_M<r_N$, are \emph{concordant}
if their rank functions satisfy $r_M(Y)-r_M(X) \leq r_N(Y)-r_N(X)$ for all $X
\subset Y \subseteq E$. 
\cite{BGW}.



\emph{Flag matroids} are an important family of \emph{Coxeter matroids} \cite{BGW}. There are several equivalent ways to define them; in particular they also have an algebro-geometric interpretation. We proceed constructively. Given pairwise concordant matroids $M_1, \ldots, M_m$ on $E$  of ranks $k_1<  \cdots <  k_m$, consider the collection of flags $(B_1, \ldots, B_m)$, where $B_i$ is a basis of $M_i$ and $B_1 \subset \cdots \subset B_m$. Such a  collection of flags is called a \emph{flag matroid}, and $M_1, \ldots, M_m$ are called the \emph{constituents} of $\F$. 


For each flag $B=(B_1, \ldots,  B_m)$ in $\F$ let $v_B=v_{B_1}+\cdots+v_{B_m}$,  where $v_{\{a_1,\dots,a_i\}}=e_{a_1}+\cdots+e_{a_i}$. The \emph{flag matroid polytope} is $P_\F= \mathrm{conv}\{v_B:B \in \F\}$.

\begin{theorem} \cite[Cor 1.13.5]{BGW} \label{th:flagpoly}
If $\F$ is a flag matroid with constituents $M_1, \ldots, M_k$,
then $P_\F = P_{M_1} + \cdots + P_{M_k}$.
\end{theorem}

As mentioned above, this implies that every flag matroid polytope is a signed
Minkowski sum of simplices $\Delta_I$; the situation is particularly nice for truncation flag matroids, 
which we now define.

Let $M$ be a matroid over the ground set $E$ with rank $r$. Define
$M_i$ to be the rank $i$ \emph{truncation} of $M$, whose bases are the independent sets of
$M$ of rank $i$. One easily checks that the truncations of a matroid are concordant,
and this motivates the following definition of Borovik, Gelfand, Vince, and White. 

\begin{deff}\cite{BGVW}
The flag $\F(M)$ with constituents $M_1,\dots,M_r$ is a flag matroid, called the \emph{truncation flag matroid} or 
\emph{underlying flag matroid} of $M$. 
\end{deff}

%







Our next goal is to present the decomposition of a truncation flag
matroid polytope as a signed Minkowski sum of simplices. For that purpose, we define the \emph{gamma invariant} of $M$ to be $\gamma(M) = b_{20}-b_{10}$, where $T_M(x,y) = \sum_{i,j} b_{ij}x^iy^j$ is the Tutte polynomial of $M$.

\begin{prop}\label{prop:gamma} The gamma invariant of a matroid is given by
\[
\gamma(M) = \sum_{I \subseteq E} (-1)^{r - |I| } {r - r(I)+1
\choose 2}.
\]
\end{prop}

\begin{proof}
We would like to isolate the coefficient of $x^2$ minus the
coefficient of $x$ in the Tutte polynomial $T_M(x,y)$.  We will
hence ignore all terms containing $y$ by evaluating $T_M(x,y)$ at
$y=0$, and then combine the desired $x$ terms through the
following operations:
\begin{eqnarray*}
\gamma(M)&:=&\frac{1}{2}\left[\frac{d^2}{dx^2}(1-x)T_M(x,0)\right]_{x=0}\\
&=&\frac{1}{2} \left[\frac{d^2}{dx^2}\sum_{I\subseteq E}(-1)^{|I|-r(I)+1}(x-1)^{r-r(I)+1}\right]_{x=0}\\
&=&\sum_{I\subseteq E}(-1)^{r-|I|}{r-r(I)+1\choose2},
\end{eqnarray*}
as we wished to show. \end{proof}

%
Unlike the beta invariant, the gamma invariant is not necessarily nonnegative. In fact its sign is not simply a function of $|E|$ and $r$.  For example, $\gamma(U_{k,n})= - {{n-3} \choose {k-1}}$, and $\gamma(U_{k,n}\oplus C)= {{n-2} \choose {k-1}}$ where $C$ denotes a coloop.

As we did with the beta invariant, define the \emph{signed gamma
invariant} of $M$ to be $\ggamma(M) = (-1)^{r(M)} \gamma(M)$.

\begin{theorem} \label{th:truncation} The truncation
flag matroid polytope of $M$ can be expressed as:
\[
P_{\F(M)} = \sum_{I\subseteq E} \ggamma(M/I)\Delta_{E-I}.
\]
\end{theorem}

\begin{proof} By Theorems \ref{th:sum} and  \ref{th:flagpoly}, $P_{\F(M)}$ is 
\[
\sum_{i=1}^r P_{M_i} 
= \sum_{i=1}^r\sum_{I\subseteq E}\sum_{J\subseteq I}(-1)^{|I|-|J|}(i-r_i(E-J))\Delta_I,
\]
where $r_i(A) = \min\{i, r(A)\}$ is the rank function of $M_i$. Then 
\begin{eqnarray*}
P_{\F(M)} & = & \sum_{I\subseteq E}\left[\sum_{J\subseteq I}\;(-1)^{|I|-|J|}\sum_{i=r(E-J)+1}^r(i-r(E-J))\right]\Delta_I \\
&=& \sum_{I\subseteq E}\left[\sum_{J\subseteq I}(-1)^{|I|-|J|}{r-r(E-J)+1\choose 2}\right]\Delta_I\\
&=& \sum_{I\subseteq E}\left[\sum_{X\subseteq I}(-1)^{|X|}{r_{M/(E-I)}-r_{M/(E-I)}(X)+1\choose 2}\right]\Delta_I\\
&=& \sum_{I\subseteq E}\ggamma(M/(E-I))\Delta_I
\end{eqnarray*}
as desired. \end{proof}

\begin{cor}\label{vol.cas.}
If a connected matroid $M$ has $n$ elements, then 
\[
\vol P_{\F(M)} = \frac{1}{(n-1)!} \sum_{(J_1,
\ldots, J_{n-1})} \ggamma(M/J_1) \ggamma(M/J_2)
\cdots \ggamma(M/J_{n-1}),
\]
summing over the ordered collections of sets
$J_1, \ldots, J_{n-1} \subseteq [n]$ such that, for any distinct $i_1,
\ldots, i_k$, $|J_{i_1} \cap \cdots \cap J_{i_k}| < n - k$.
\end{cor}

\begin{proof} This follows from Proposition \ref{prop:gen.mixed.vol}  and Theorem \ref{th:truncation}.
\end{proof}

\begin{exa}
Let $M$ be the matroid on $[3]$ with bases $\{1,2\}$ and $\{1,3\}$.
The flags in $\F(M)$ are:
$
\{1\}\subseteq\{1,2\},
\{1\}\subseteq\{1,3\},
\{2\}\subseteq\{1,2\}, 
\{3\}\subseteq\{1,3\},
$
so the vertices of $P_{\F(M)}$ are $(2,1,0), (2,0,1), (1,2,0),(1,0,2)$, respectively. Theorem \ref{th:truncation} gives
$P_{\F(M)}=\Delta _{123}+\Delta _{23}. $
Since $\widetilde{\gamma}(M) = \widetilde{\gamma}(M/1)=1$, Corollary \ref{vol.cas.} gives
\[
\vol P_{\F(M)}=\dfrac{1}{2!}[\widetilde{\gamma} (M/
\emptyset)\widetilde{\gamma} (M/ \emptyset)+\widetilde{\gamma} (M/
\emptyset)\widetilde{\gamma} (M/ 1)+\widetilde{\gamma} (M/
1)\widetilde{\gamma} (M/ \emptyset)]=\dfrac{3}{2}.
\]
\end{exa}

\section{\textsf{acknowledgments.}}
This work is part of the SFSU-Colombia Combinatorics Initiative. We are very grateful to SFSU and the Universidad de Los Andes for supporting this initiative.

\footnotesize{

}

\end{document}